\documentclass[a4paper, 12pt, reqno]{amsart}
\usepackage{amsmath, amssymb, latexsym, graphicx, hyperref}

\newcommand{\N}{\mathbb{N}}
\newcommand{\R}{\mathbb{R}}

\begin{document}

\title{Correcting a minor error in Cantor's \\
calculation of the power of the continuum}

\maketitle

\begin{center}
\vspace{-18pt} \small{Antonio Leon
(\href{mailto:aleons@educa.jcyl.es}{aleons@educa.jcyl.es})\\I.E.S.
Francisco Salinas. Salamanca, Spain.}
\end{center}

\pagestyle{myheadings}

\markboth{\small{Correcting a minor error in Cantor's calculation
of the power of the continuum}}{\small{Correcting a minor error in
Cantor's calculation of the power of the continuum}}

\begin{abstract}
Cantor's algebraic calculation of the power of the continuum
contains an easily repairable error related to Cantor own way of
defining the addition of cardinal numbers. The appropriate
correction is suggested.
\end{abstract}

\section{The exponentiation of powers}

\noindent Cantor's most significant contribution to the theory of
transfinite numbers is, without a doubt, \emph{Beitr\"age zur
Begr\"undung der transfiniten Mengelehre}\footnote{Translated to
English by P. E. B. Jourdain in 1915 as Contributions to the
Founding of the Theory of Transfinite Numbers, \cite{Cantor1955}}.
A memory of more than 70 pages divided into two parts which
appeared in the Mathematische Annalen in the years 1895 and 1897
respectively (\cite{Cantor1895}, \cite{Cantor1897}).
\emph{Beitr\"age}'s first six epigraphs are devoted to found the
arithmetics of cardinals. Cantor begins by defining the concept of
set and the union of disjoint sets, after which he proposes the
following definition of power or cardinal number
(\cite{Cantor1955}, p. 86):\\

\begin{quote}
\small{

We call by the name ''power" or ''cardinal number" of [the set] M
the general concept which, by means of our active faculty of
thought, arises from the set $M$ when we make abstraction of the
nature of its various elements $m$ and of the order in which they
are give.

\indent \\We denote the result of this double act of abstraction,
the cardinal number or power of $M$ by $\overline{\overline{M}}$.

}
\end{quote}

\noindent \\Cantor continues by defining the concept of
equivalence for sets: two sets $M$ and $N$ are said equivalent,
symbolically $M \sim N$, if they can be put into a one to one
correspondence (\cite{Cantor1895}, p. 86). He then proves that two
sets are equivalent if, and only if, they have the same power
(\cite{Cantor1895}, pp. 87-88). After extending the notions of
''greater than" and ''less than" to cardinals numbers, Cantor
defines in set theoretical terms the operations of addition and
multiplication of cardinals (\cite{Cantor1895}, pp. 91-94). Since
multiplication cannot be easily extended to the case of infinitely
many factors, Cantor defines the notion of covering in order to
define the exponentiation of (finite and transfinite) powers
(\cite{Cantor1955}, p. 94):\\

\begin{quote}

\small{

By a ''covering of the set $N$ with elements of the set M," or
more simply, by a ''covering of $N$ with M," we understand a law
by which with every element $n$ of $N$ a definite element of $M$
is bound up, where one and the same element of $M$ can come
repeatedly into application. The element of $M$ bound up with $n$
is, in a way, a one value function of $n$, and may be denoted by
$f(n)$; it is called a ''covering function of $n$." The
corresponding covering of $N$ will be called $f(N)$.

}

\end{quote}

\noindent \\So, if $N = \{a, b, c, d, e\}$ and $M=\{0,1\}$ the
coverings of $N$ are:

\begin{equation}
f(N) = 10110, \ f'(N) = 00111, \ f''(N) = 10111 \dots
\end{equation}

\noindent \\The totality of different coverings of $N$ with M,
denoted as $(N | \: M)$, forms a set Cantor called the
''covering-set of $N$ with $M$":

\begin{equation}
(N | \: M) = \{f(N), f'(N), f''(N), \dots \}
\end{equation}

\noindent For instance:

\begin{equation}
(\{a, b, c\} | \: \{0,1\}) = \{000, 001, 010, 100, 011, 101, 110,
111\}
\end{equation}

\noindent \\It is immediate that, if $N'$ and $M'$ are two sets
such that $N \sim N'$ and $M \sim M'$ it holds $(N | \: M) \sim
(N' | \: M')$. This equivalence is necessary to prove the
cardinality of $(N | \: M)$ depends exclusively upon the cardinal
$a$ of $M$ and $b$ of $N$, which in turn makes it possible to
define the cardinal number $a^b$ as the cardinal of the set $(N |
\: M)$ (\cite{Cantor1895}, p. 95):

\begin{equation}\label{eqn:def a^b}
a^b =(\overline{\overline{N |\: M}})
\end{equation}

\noindent \\It immediately follows from the above definition that
(\cite{Cantor1895}, p. 95):

\begin{align}
\left( (N |\: M) \cdot (P |\:M)\right) \sim \left((N,P) |\:
M)\right)\\
\left( (P |\:M) \cdot (P |\:N)\right) \sim \left(P |\:
(M \cdot N)\right)\\
\left( P |\: (N |\: M)\right) \sim \left( (P \cdot N) |\: M
\right)
\end{align}

\noindent \\Thus, if $\overline{\overline{N}} =a$,
$\overline{\overline{M}} = b$, and $\overline{\overline{P}} = c$,
we will have (\cite{Cantor1895}, p. 95):

\begin{align}
a^b \cdot a^c &= a^{b+c}\\
a^c \cdot b^c &= (a \cdot b)^c\\
(a^b)^c &= a^{b \cdot c}
\end{align}

\noindent \\The main consequence of the concept of covering is,
therefore, the possibility of defining the exponentiation of
cardinals numbers even for infinite values of the exponent.\\

Particularly important in Cantor's determination of the power of
the continuum is the covering-set of the set $\N$ of all finite
cardinals $\{ 1$, 2, 3, $\ldots \}$ with the set $M$ of two
elements $\{0,1\}$. This covering-set is the set of all binary
infinite strings of 1s and 0s:

\begin{align*}
&1110001110001110001101011 \ldots\\
&1000000000000011111110000 \ldots\\
&0001110001110001111111111 \ldots\\
&\vdots
\end{align*}

\noindent Each of these strings, when preceded by a decimal point,
is the binary representation of a real number $x$ in the real
interval $[0,1]$.\\[14pt]

\section{The power of the continuum}

\noindent In the year 1872 Cantor published a paper on the
continuity and infiniteness of the set $\R$ of real numbers
\cite{Cantor1874}. In this paper, Cantor considered axiomatic the
one to one correspondence between the real line points (linear
continuum) and the real numbers. Two years after, in 1874, he gave
his first proof on the non enumerable nature of the set $\R$ of
real numbers \cite{Cantor1874} (his second proof -the diagonal
method- on the existence of non denumerable sets was published in
1891 \cite{Cantor1891}). In 1878 he proved the equivalence between
linear and n-dimensional continuums \cite{Cantor1878}. All these
results were directly or indirectly necessary to prove the power
of the continuum is the transfinite cardinal $2^{\aleph_0}$
(\cite{Cantor1955}, p. 96):\\

\begin{quote}

\small{

We see how pregnant and far-reaching these simple formul{\ae}
extend to powers are by the following example. If we denote the
power of the linear continuum X (that is the totality X of real
numbers $x$ such that $x \geq$ and $\leq 1$) by $c$, we easily see
that it may be represented by, among others, the formula:

\begin{equation}
c = 2^{\aleph_0}
\end{equation}

\noindent where \S 6 gives the meaning\footnote{In \S 6 Cantor
defines $\aleph_0$ as the smallest transfinite cardinal number:
the cardinal of the set of all finite cardinals: $\aleph_0 =
\overline{\overline{\{\nu\}}}$} of $2^{\aleph_0}.$. In fact, by
(\ref{eqn:def a^b}), $2^{\aleph_0}$ is the power of all
representations

\begin{equation}
x = \frac{f(1)}{2} + \frac{f(2)}{2^2} + \dots +
\frac{f(\nu)}{2^\nu} + \dots \text{( where $f(\nu)$ = 0 or 1)}
\end{equation}

\noindent of the numbers $x$ in the binary system. If we pay
attention to the fact that every number $x$ is only represented
once, with the exception of the numbers $x = \frac{2\nu+1}{2^\mu}
< 1$, which are represented twice over, we have, if we denote the
''enumerable" totality of the latter by $\{s_\nu\}$,

\begin{equation}\label{eqn:(sv, X)}
2^{\aleph_0} = (\overline{\overline{\{s_\nu\}, \:X }})
\end{equation}

\noindent If we take away from $X$ any ''enumerable" set
$\{t_\nu\}$ and denote the remainder by $X_1$, we have:

\begin{equation}
X = (\{t_\nu\}, X_1) = (\{t_{2\nu-1}\}, \{t_{2\nu}\}, X_1)
\end{equation}

\begin{equation}
(\{s_\nu\}, X) = (\{s_\nu\}, \{t_\nu \}, X_1)
\end{equation}

\begin{equation}
\{t_{2\nu-1}\} \sim \{s_\nu\}, \ \{t_{2\nu}\} \sim \{s_\nu\}, \
X_1 \sim X_1
\end{equation}

\noindent so

\begin{equation}
X \sim (\{s_\nu\}, X)
\end{equation}

\noindent and thus

\begin{equation}
2^{\aleph_0} = \overline{\overline{X}} = c\\[20pt]
\end{equation}

}

\end{quote}

\noindent \\According to Cantor's notation, $(M, N)$ is the union
of two sets $M$ and $N$ which have no common elements
(\cite{Cantor1895}, p. 85). The union of disjoint sets is
essential in Cantor's definitions of arithmetic operations. So, if
$\overline{\overline{M}} =a$ and $\overline{\overline{N}} = b$,
the sum of the cardinals $a$ and $b$ is given by
(\cite{Cantor1895}, p. 91):

\begin{equation}
a + b = (\overline{\overline{M, N}})
\end{equation}

\noindent \\being $M \cap N = \emptyset$. It is therefore clear
the meaning of the above equation (\ref{eqn:(sv, X)})

\begin{equation}
2^{\aleph_0} = (\overline{\overline{\{s_\nu\}, \:X }})
\end{equation}

\noindent \\Although being $\{s_\nu\}$ the enumerable totality of
real numbers $x = \frac{2\nu+1}{2^\mu} < 1$, which have two binary
representations (\cite{Cantor1895}, p. 96), and $X$ the set of of
all real numbers in $[0, \: 1]$ (\cite{Cantor1895}, p. 96), it is
also clear that $\{s_\nu\} \subset X$, so that $\{s_\nu\} \cap X
\neq \emptyset$. This little difficulty is easily solved by
redefining the sets involved in the proof. In fact, let $B$ be the
set of all binary infinite strings of 0s and 1s (the covering-set
of $\N$ by $\{0, \: 1\}]$). By definition, we have

\begin{equation}
\overline{\overline{B}}= 2^{\aleph_0}
\end{equation}

\noindent  \\The set $B$ can be divided into two disjoint sets:
the set $B_X$ and the set $B_S$. The set $B_X$ is the set of all
binary strings representing all real numbers of the set $X = [0,
\: 1]$ except the strings of the second binary expressions of all
$X$'s elements which have two binary expressions (the first one
ending by an infinite string of 0s, and the second by an infinite
string of 1s). The set $B_S$ is just the denumerable set of all
those second binary strings. Evidently, we will have:

\begin{equation}
B = B_X \cup B_S
\end{equation}

\noindent and:

\begin{equation}
B_X \sim X \sim \R
\end{equation}

\noindent \\being $B_X$, $X$ and $\R$ non denumerable. From this
point we only have to follow Cantor's argument. Let $T$ be any
denumerable subset of $B_X$ and let $B'_X$ be the complement of
$T$ with respect to $B_X$, i.e. $B'_X = B_X - T$, we can write:

\begin{equation}\label{eqn:Bx = T U B'x}
B_X = T \cup B'_X
\end{equation}

\noindent \\being $T \cap B'_X = \emptyset$. Since $T$ is
denumerable its elements can be indexed by the totality of natural
numbers. Consequently we can consider two disjoint denumerable
subsets $T_E$ and $T_O$, whose elements are respectively indexed
by the even and the odd natural numbers. Equation (\ref{eqn:Bx = T
U B'x}) can then be rewritten as:

\begin{equation}\label{eqn:Bx = TE U TO U B'X}
B_X = T_E \cup T_O \cup B'_X
\end{equation}

\noindent \\where $T_E \cup T_O = T$; $T_E \cap T_O = \emptyset$.
From (\ref{eqn:Bx = T U B'x}) we also get:

\begin{equation}\label{BS U BX = BS U T U B'X}
B_S \cup B_X = B_S \cup T \cup B'_X
\end{equation}

\noindent and being

\begin{align}
T_E &\sim B_S &&\text{($T_E$ and $B_S$ are denumerable)}\\
T_O &\sim T &&\text{($T_O$ and $T$ are denumerable)}\\
B'_X &\sim B'_X &&\text{(Every set is equivalent to itself)}\\
\end{align}

\noindent we have:

\begin{equation}
T_E \cup T_O \cup B'_X \sim B_S \cup T \cup B'_X
\end{equation}

\noindent \\and then, in accordance with (\ref{eqn:Bx = TE U TO U
B'X}) and (\ref{BS U BX = BS U T U B'X}):

\begin{equation}
B_X \sim B_S \cup B_X
\end{equation}

\noindent \\and taking into account that $B_S \cup B_X = B$, we
get:

\begin{equation}
B_X \sim B
\end{equation}

\noindent and then

\begin{equation}
\overline{\overline{B}}_X = \overline{\overline{B}} = 2^{\aleph_0}
\end{equation}

\noindent \\Finally, being $B_X \sim X \sim \R$, we can write:

\begin{equation}
\overline{\overline{X}} = \overline{\overline{\R}} = 2^{\aleph_0}
\end{equation}

\noindent \\which proves that, in fact, the power of the continuum
is the transfinite cardinal $2^{\aleph_0}$.


\begin{thebibliography}{1}

\bibitem{Cantor1874}
Georg Cantor, \emph{{\"{U}}ber eine eigenschaft aller reallen
algebraishen
  zahlen}, Journal f{\"{u}}r die reine und angewandte Mathematik \textbf{77}
  (1874), 258--262.

\bibitem{Cantor1878}
\bysame, \emph{Ein {B}eitrag zur {M}annigfaltigkeitslehre},
Journal f{\"{u}}r
  die reine und angewandte Mathematik \textbf{84} (1878), 242--258.

\bibitem{Cantor1891}
\bysame, \emph{{\"{U}}ber {E}ine elementare frage der
mannigfaltigkeitslehre},
  Jahresberich der {D}eutschen {M}athematiker {V}ereiningung, vol.~1, 1891.

\bibitem{Cantor1895}
\bysame, \emph{Beitr{\"{a}}ge zur {B}egr{\"{u}}ndung der
transfiniten
  {M}engenlehre}, Mathematische Annalen \textbf{XLVI} (1895), 481 -- 512.

\bibitem{Cantor1897}
\bysame, \emph{Beitr{\"{a}}ge zur {B}egr{\"{u}}ndung der
transfiniten
  {M}engenlehre}, Mathematishe Annalen \textbf{XLIX} (1897), 207 -- 246.

\bibitem{Cantor1955}
\bysame, \emph{Contributions to the founding of the theory of
transfinite
  numbers}, Dover, New York, 1955.

\end{thebibliography}

\providecommand{\bysame}{\leavevmode\hbox
to3em{\hrulefill}\thinspace}
\providecommand{\MR}{\relax\ifhmode\unskip\space\fi MR }
\providecommand{\MRhref}[2]{%
  \href{http://www.ams.org/mathscinet-getitem?mr=#1}{#2}
} \providecommand{\href}[2]{#2}

\end{document}